%%%%%%%%%%%%%%%%%%%%%%%%%%%%%%%%%%%%%%%%%%%%%%%%    
%
%        THIS IS A  PLAIN TeX FILE
%
%%%%%%%%%%%%%%%%%%%%%%%%%%%%%%%%%%%%%%%%%%%%%%%%

\magnification=1200

%\font\AAA=Times.dfont  at 12pt
 %\font\BBB=Times.dfont at 8pt

%\font\AAA=cmr10 at 12pt
%\font\BBB=cmr10 at 8pt

\overfullrule=0in

\def\boxit#1{\hbox{\vrule
 \vtop{%
  \vbox{\hrule\kern 2pt %
     \hbox{\kern 2pt #1\kern 2pt}}%
   \kern 2pt \hrule }%
  \vrule}}

  \def\harr#1#2{\ \smash{\mathop{\hbox to .3in{\rightarrowfill}}\limits^{\scriptstyle#1}_{\scriptstyle#2}}\ }

  \def\leftharr#1#2{\ \smash{\mathop{\hbox to .3in{\leftarrowfill}}\limits^{\scriptstyle#1}_{\scriptstyle#2}}\ }

\def\ss{\subset}

\def\dim{{\rm dim}}

\def\max{{\rm max}}

\def\arr{\longrightarrow}
\def\supp{{\rm supp}}
\def\Link{{\rm Link}}

\def\Theorem#1{\medskip\noindent {\bf THEOREM \bf #1.}}
\def\Prop#1{\medskip\noindent {\bf Proposition #1.}}
\def\Cor#1{\medskip\noindent {\bf Corollary #1.}}
\def\Lemma#1{\medskip\noindent {\bf Lemma #1.}}
\def\Remark#1{\medskip\noindent {\bf Remark #1.}}
\def\Note#1{\medskip\noindent {\bf Note #1.}}
\def\Def#1{\medskip\noindent {\bf Definition #1.}}

\def\Ex#1{\medskip\noindent {\bf Example \bf    #1.}}
\def\Qu#1{\medskip\noindent {\bf Question \bf    #1.}}

\def\pf{\medskip\noindent {\bf Proof.}\ }
\def\qed{\hfill  $\vrule width5pt height5pt depth0pt$}

   \def\cp{{\cal P}}   
   \def\co{{\cal O}}
\def\ce{{\cal E}}

\def\cd{{\cal D}}

\def\cp{{\cal P}}

\def\vf{\varphi}

\def\wt{\widetilde}

\def\and{\qquad {\rm and} \qquad}
\def\arr{\longrightarrow}
\def\ol{\overline}
\def\bbr{{\bf R}}
\def\bbc{{\bf C}}
\def\bbq{{\bf Q}}
\def\bbz{{\bf Z}}
\def\bbp{{\bf P}}

\def\a{\alpha}
\def\b{\beta}

\def\o{\omega}

\def\s{\sigma}
\def\x{\xi}

\def\L{\Lambda}
\def\G{\Gamma}

\def\dbar{\ol{\partial}}

\def\CC{3}
\def\DD{5}
\def\EE{4}

\font\titfont=cmr10 at 12pt

\def\Scar{{\rm \Sigma_{\G}}}

\centerline{\titfont BOUNDARIES OF POSITIVE HOLOMORPHIC CHAINS}\smallskip

\centerline{\titfont AND THE RELATIVE HODGE QUESTION}

\medskip

\centerline{\bf by}

\medskip

\centerline{\titfont F. Reese Harvey and H. Blaine Lawson, Jr.*}

 \vskip .3in

\smallbreak\footnote{}{ $ {} \sp{ *}{\rm Partially}$  supported by
the N.S.F. }

\centerline{\bf Abstract} \medskip
  \font\abstractfont=cmr10 at 10 pt

{{\parindent=.8in\narrower\abstractfont  \noindent
We characterize the boundaries of positive holomorphic chains in an
arbitrary complex manifold.
\smallskip
\noindent
We then consider a compact oriented real submanifold  of dimension
$2p-1$ in a compact K\"ahler manifold $X$ and address the question
of which relative homology classes in $H_{2p}(X,M;\,\bbz)$ are 
represented by positive holomorphic chains. Specifically,
we define what it means for a class $\tau\in H_{2p}(X,M;\,\bbz)$
to be of type (p,p) and positive. It is then shown that $\tau$ 
has these properties if and only if $\tau = [T+S]$ where $T$ 
is a positive holomorphic chain with $dT=\partial \tau$ and 
$S$ is a positive (p,p)-current with $dS=0$.

}}
 
 \vskip.3in

\def\om{{1\over p!}\omega^p} 
 \def\Ga{{\G}}
\def\bbt{{ V}}

\noindent
{\bf \S 1. Introduction.}  
In the first part of this note we establish  a general result concerning boundaries of
 positive holomorphic chains in a complex manifold $X$.   
 In the second part we address the ``Relative Hodge Question'': {\sl  When is a homology
 class $\tau\in H_{2p}(X,M;\,\bbz)$ represented by a positive holomorphic chain?}
 Assuming $M$ is a real $(2p-1)$-dimensional submanifold we are able to give a surprisingly
 full answer.
 
 We begin our discussion of the first part by presenting  some interesting special cases  which are quite different in nature.  The first main theorem is then formulated and proved in  Section 2.

 To start, suppose  $X$   compact   and let  
 $\Ga$  be a current of dimension $2p-1$ in $X$.  By a positive holomorphic
 $p$-chain with boundary $\Ga$ we mean a finite sum $\bbt=\sum_k m_k V_k$
with $m_k\in \bbz^+$ and $V_k$ an irreducible complex analytic variety of dimension $p$ and finite volume in $X-\supp\Ga$, such that $d\bbt=\Ga$ as currents on $X$.

Equip $X$ with a hermitian metric and let $\o$ denote its associated $(1,1)$-form.
A    real  $(2p-1)$-form $\a$  will be called
a {\bf $(p,p)$-positive linking form} if 
$$
\qquad\qquad \qquad   d^{p,p}\a+\om\geq 0 \qquad\qquad {\rm  (strongly\  positive) }
$$
where  $ d^{p,p}\a$ denotes the $(p,p)$-component of $d\a$. 
(See [HK] or [H$_2$] for the  definition  of   strongly and  weakly positive currents.)
The numbers   $\int_\G \a$ with $\a$ as above, will be called the {\bf $(p,p)$-linking numbers} of $\G$.

\eject
\Theorem{1.1} {\sl Let $\Ga = \sum_{k=1}^N n_k \Ga_k$ be an integer linear combination of compact, mutually disjoint, $C^1$-submanifolds of dimension $2p-1$  in $X$, each of which has a real analytic point.
Then $\Ga = d\bbt$ where $\bbt$ is a positive holomorphic p-chain   if and only if  the 
$(p,p)$-linking numbers of $\G$ are bounded below.
} 
 
 \medskip

  \Note{1.2}   The condition that the linking numbers of $\G$ are bounded below is easily seen to be 
 independent of the choice of hermitian metric on $X$.     However, for any given metric we have the  precise statement that  $\G$ bounds a  positive holomorphic chain of mass $\leq \Lambda$ if and only if
$$
\int_\G \a\ \geq\ -\Lambda \qquad {\rm for\ all \ } (p,p){\rm -positive\ linking\ forms\ }\ \a
%\ {\rm with\ compact\ support.} 
\eqno{(1.1)} 
$$

\Note{1.3}  We shall actually prove the theorem in the more general situation where
$\Ga$ is allowed to have a  ``scar"  set and the real analyticity assumption is replaced
by  a weaker ``push-out'' hypothesis (see section 2).  When $p>1$, this hypothesis is satisfied at any point where the boundary is smooth and its Levi form has at least one negative eigenvalue.
In all these cases, one has  regularity at almost all  points of $\Ga$.  This boundary regularity  is discussed in  [HL$_1$] and [H$_2$].

\Remark{1.4}  When $X$ is a projective surface and $p=1$, a much stronger result is
conjectured: namely, $\Ga$ bounds a positive holomorphic 1-chain if and only if
$$
\int_{\Ga}  d^c u\ \geq -\Lambda \qquad {\rm for\  all\ } u\in C^\infty(X) { \rm \ with\ }  dd^cu+\o \geq 0.
\eqno{(1.2)}
$$
Functions $u$  with  $dd^cu+\o \geq 0$ are called {\sl quasi-plurisubharmonic}.  They were introduced
by Demailly and play an important role in complex analysis [D], [GZ]. Condition (1.2) is equivalent to the condition that
$$
{1\over\ell} \Link_{\bbp}(\Ga, Z)\ \geq \ -\L \qquad {\rm for\  all\ positive\  divisors\ }  Z \ {\rm in\ } X-\Ga
$$
of sections $\s\in H^0(X,\co{}(\ell))$, $\ell >0$, where $\Link_{\bbp}$ denotes the  {\sl  projective linking number } introduced in [HL$_5$].
In this form the conjecture extends to all dimensions and codimensions (for $X$ projective) and is a consequence of the above case:  $p=1$ in surfaces. 
%This represents a projective version of a theorem of Alexander and Wermer [AW].  
All this is established in [HL$_5$] where the conjectures are also related to the projective hull introduced in [HL$_4$].

Although the hypothesis of Theorem 1.1  is conjecturally  too strong  for projective manifolds,  it does give the ``correct'' result in the  general case.  For example, if $X$ is a non-algebraic 
K3-surface, there appears to be no simpler condition characterizing the boundaries
of positive holomorphic 1-chains.

\Remark{1.5}  The   Linking Condition (1.1) forces the components of $\Ga$ to be maximally complex CR-manifolds.   Maximal complexity is equivalent to the assertion that $\Ga = \Ga_{p-1,p} + \Ga_{p,p-1}$  where $ \Ga_{r,s} $ denotes the Dolbeault component
of $\Ga$ in bidimension $(r,s)$. To see that this must hold,  note that any $\a\in\ce^{r,2p-1-r}(X)$ with $r\neq p-1,p$
satisfies $d^{p,p}\a+\o \geq 0$ since $d^{p,p}\a=0$.

\medskip

 Theorem 1.1 extends to characterize boundaries of compactly supported holomorphic chains in certain non-compact spaces.  A complex $n$-manifold $X$ is called  {\sl $q$-convex} if there exists a proper exhaustion function $f:X\to \bbr^+$ such that $dd^c f$ has
at least $n-q+1$ strictly positive eigenvalues outside some compact subset of $X$.

\Theorem{1.6}   {\sl  Theorem 1.1 remains valid (for compactly supported holomorphic chains  $V$) in any $q$-convex hermitian manifold with
$q\leq p$.  }

\Remark{1.7}  If $X$ is 1-convex (i.e., strongly pseudoconvex), then Theorem 1.1 is valid for all $p$.
If, further, $X$ admits a proper exhaustion which is strictly plurisubharmonic  everywhere (i.e., $X$ is Stein), much stronger results are known.  Condition (1.1) implies maximal complexity, and for
$p>1$ this condition alone implies that $\Ga$ bounds a holomorphic $p$-chain [HL$_1$].  Condition
(1.1) also implies the {\sl moment condition}: $\Ga (\a) =0$ for all $(p,p-1)$-forms $\a$ with 
$\dbar \a=0$. When $p=1$ this implies that $\Ga$ bounds a holomorphic $1$-chain [HL$_1$]. 
 
 Analogous remarks apply to results of [HL$_2$] in the $q$-convex spaces $\bbp^n-\bbp^{n-q}$.

\Remark{1.8}   Condition (1.1) implies that $\int_{\Ga}  \a\ \geq 0$ for all $\a$ with $d^{p,p}\a \geq 0$.
If $X$ is a Stein manifold embedded in some $\bbc^n$, this in turn implies that 
the linking number Link$(\Ga, Z)\geq0$ for all algebraic subvarieties $Z$ of codimension\ $p$
in $\bbc^n-\Ga$. By Alexander-Wermer [AW], [W$_2$] this last condition alone implies that
$\Ga$ bounds a positive holomorphic $p$-chain in $X$.

Theorem 1.1 also holds ``locally'', that is, it extends to any non-compact hermitian manifold
$X$ where neither $\G$ nor $V$ are assumed to have compact support.

\Theorem{1.9}  {\sl  Suppose $X$ is a non-compact hermitian manifold, and let $\G=\sum_j n_j \G_j$
be a locally finite integral combination of disjointly embedded $C^1$-submanifolds of dimension
$2p-1$, each of which has a real analytic point.  Then $\G$ is the boundary of a holomorphic
$p$-chain $V$ of mass $M(V) \leq \Lambda$ (whose support is a closed but not necessarily compact
analytic subvariety  of $X-\supp \G$)  
if and only if $\int_\G \a \geq -\Lambda$ for all $(p,p)$-positive linking forms $\a$ with compact support on $X$.
}

\medskip

In the last section  of this paper we further  weaken our  hypotheses on $\G$
to an assumption that each component $\G_k$ be {\sl residual} at some point.
(See \S 3 for the definition.)  The concept of residual submanifolds leads to questions
of some independent interest.

\medskip
In Section three we address a  question related to the Characterization Theorems
above. Let $j:M\ss X$ be a compact oriented real submanifold of dimension $2p-1$
in a compact K\"ahler manifold $X$.  Represent the relative homology group
$H_{2p}(X,M;\,\bbr)$ by $2p$-currents $T$ on $X$ with $dT=j_*S$ for some $(2p-1)$-current $S$  
on $M$. One can ask: {\sl When does a given class 
$\tau\in H_{2p}(X,M;\,\bbr)$ contain a positive holomorphic chain?}

As a first step we show that for every $T$ as above and every $d$-closed form
$\vf$ on $X$ the pairing $T(\vf)$ depends only on the relative class $\tau=[T]$.
This allows us to introduce a real Hodge filtration  on $H_{2p}(X,M;\,\bbz)_{\rm mod\ tor}$
which extends the standard one on the subgroup $H_{2p}(X;\,\bbz)_{\rm mod\ tor}$.
It also allows us to formulate the following.

\Def{1.10}  A class $\tau\in H_{2p}(X,M;\,\bbr)$ is a {\bf positive (p,p)-class} if 
$\tau (\vf)\geq0$ for all $2p$-forms $\vf$ with $d\vf=0$ and $\vf^{p,p}\geq 0$.

\vfill\eject

\Theorem {1.11}  {\sl Let $M\ss X$ be as above and suppose each
component of $M$ has a real analytic point. 
Let $\tau \in   {H}_{2p}(X, M;\,\bbz)_{\rm mod\ tor}$ be a positive $(p,p)$-class.
Then there exists a positive holomorphic $p$-chain $T$ on $X$
with $d T=\partial \tau$ and a positive
$(p,p)$-current $S$ with $dS=0$ such that $\tau =[T+S]$.

In particular, if the positive  classes in $H_{p,p}(X;\,\bbq)$ are  
represented  by positive holomorphic chains with rational coefficients, then so are all the positive
 classes in $H_{p,p}(X, M;\,\bbq)$.}

\Remark{1.12} This last result is a strengthening of the previous ones (in the K\"ahler case).
Let $\tau$ be as in Theorem 1.11 and note that  $\G=\partial \tau =\sum_k n_k [M_k]$ 
where $M_1,...,M_\ell$ are the connected components of $M$ and the $n_k$'s are integers.
If $\tau$ is a positive $(p,p)$-class, then  $\tau(d\a +{1\over p!}\o^p)\geq0$ 
 whenever $d^{p,p}\a +{1\over p!}\o^p\geq0$.  Therefore for any $(p,p)$-positive linking form
 $\a$ we have $\G(\a) = (\partial \tau)(\a) = \tau(d\a)=\tau(d^{p,p}\a)= 
\tau(d^{p,p}\a+{1\over p!}\o^p) - \tau({1\over p!}\o^p)\geq -\tau({1\over p!}\o^p)$, 
and we conclude from Theorem 1.1 that $\G$ bounds a positive holomorphic 
$p$-chain $T$.  Theorem 1.11 asserts that, moreover, the absolute class $\tau-[T]$
is represented by  a positive $(p,p)$-current.

\vskip .3in
\noindent
{\bf \S 2. The   Characterization  Theorem.}  
In this section we  prove a general  theorem  which implies all of the results discussed in \S 1 except Theorem 1.11. We shall assume throughout that  $X$ is a  hermitian manifold which is not necessarily compact. 

\Def{2.1} Suppose there exists a closed subset $\Scar$ of Hausdorff $(2p-1)$-measure zero and an oriented, properly embedded, $(2p-1)$-dimensional  $C^1$ submanifold   of  
 $X-\Scar $ with connected components $\G_1,\G_2,...$.  If for given integers $n_1, n_2,...$, 
 $$
 \Ga = \sum_{k=1}^{\infty} n_k \Ga_k
 $$
 defines a current  of locally finite mass in$X$ which is $d$-closed, then $\G$
will be called a  {\bf   scarred $(2p-1)$--cycle of class $C^1$} in $X$. 
By a unique choice of orientation on   $\Ga_k$ we  assume each  $n_k > 0$.

\Ex{} Any  real analytic $(2p-1)$-cycle  is automatically a scarred $(2p-1)$-cycle (see [F, p. 433]).

\Def{2.2}  By a {\bf positive holomorphic $p$-chain with boundary $\G$} in $X$ we mean a   sum $V=\sum_k m_k V_k$ with $m_k\in \bbz^+$ and $V_k$ an irreducible $p$-dimensional complex analytic subvariety of $X-\supp \G$ such that $V$ has locally finite mass  in $X$ and $dV=\G$ as currents.

\Def{2.3} Suppose $\G$ is an embedded $(2p-1)$-dimensional oriented submanifold of a complex 
manifold. We say that $\G$ can be {\bf pushed out at } $p\in \G$ if there exists a complex
$p$-dimensional submanifold-with-boundary $(V, -\G)$ containing the point $p$ (i.e., 
$\partial V=-\G$ as oriented manifolds). 

\medskip
Our main result is the following.
\eject

\Theorem{2.4} {\sl Let $\Ga$ be a  scarred $(2p-1)$-cycle of class $C^1$ 
in  $X$ such that each component $\Ga_k$  can be pushed out at some point. Then $\Ga = d\bbt$ where $\bbt$ is a positive holomorphic p-chain with mass $M(V) \leq \Lambda$   if and only if  the 
$(p,p)$-linking numbers of $\G$ are bounded below by $-\Lambda$.}

\Remark{2.5}  We say $\G$ is {\bf two sided at $p$} if there exists a 
complex $p$-dimensional submanifold $V$ near $p$ with $\G\subset V$ near $p$.
Note that if $\G$ is real analytic and maximally complex  at $p$, then $\G$ is two-sided at $p$. 
Note also that if $\G$ is two-sided
at $p$, then $\G$ can be pushed out at $p$.\medskip

The proof of Theorem 2.4 has two parts.  First the linking condition is shown to be equivalent
to the existence of a weakly positive current $T$ of bidimension p,p satisfying $dT=\G$.
In the second part it is shown that the existence of a positive   $T$ with $dT=\G$ together with
the pushout hypothesis on $\G$ implies the existence of a positive 
holomorphic chain with boundary $\G$.

\bigskip
\centerline{\bf Solving $dT=\G$ for $T$ positive.}

\Theorem{2.6} {\sl Let $\G \in \cd_{2p-1}'(X)$ be an arbitrary current of dimension $2p-1$ on $X$.
Then $dT=\G$ for some weakly positive $(p,p)$-current with mass $M(T) \leq \Lambda$
if and only if the linking condition
$$
\int_{\G} a\ \geq \ -\Lambda
\eqno{(2.1)}
$$
is satisfied for all compactly supported, strongly positive $(p,p)$-linking forms $\a$ on $X$.
}

\pf
Let

\centerline{
$
S \ \equiv\ \{\a\in \cd^{2p-1}(X) : d^{p,p}\a+{1\over p!} \o^p\geq 0\ \   ({\rm strongly \ positive})\}
$
}
\noindent
and let
$$
C\ \equiv \ \{\G\in  \cd_{2p-1}'(X) : \G=dT\ \ {\rm for\ some\ } \ T\geq 0\ \ {\rm (weakly \ positive)
\ with\ } M(T)\leq 1\}
$$
It suffices to prove the theorem for $\Lambda =1$.  In this case the theorem states that $\G\in C$ if and only if $\G\in S^0$, where $S^0\equiv \{\G\in  \cd_{2p-1}'(X)  : \G(\a)\geq -1 $ for all $\a\in  S\}$ is the {\sl polar} of $S$.   So we must prove that
$$
C\ =\ S^0.
$$
Note that $C$ is a closed convex set in  $\cd_{2p-1}'(X)$ since the set of weakly positive $(p,p)$-currents $T$ with $M(T)\leq 1$ is compact in $ \cd_{p,p}'(X)$.  Hence by  the Bipolar Theorem [S] $C=(C^0)^0$, and it will suffice to prove that $C^0=S$.

To see this first note that
$$
T\left( d^{p,p}\a+{1\over p!} \o^p\right)\ =\ (dT)(\a)+M(T)
\eqno{(2.2)}
$$
for all weakly positive $(p,p)$-currents $T$ and all $\a\in \cd^{2p-1}(X)$.
If, in addition, $\a\in S$ and $\G\in C$, then $0\leq \G(\a)+M(T)=\G(\a)+1$, so that $S\subseteq C^0$.

It remains to show that $C^0\subseteq S$. Choose $\G=dT$ with $T=\delta_x \x$ where
$\x$ is a weakly positive $(p,p)$-vector of mass norm one at $x\in X$. Note that $\G\in C$.  By (2.2)
we have $( d^{p,p}\a+{1\over p!} \o^p)_x(\x) = \G(\a)+M(T) = \G(\a)+1$. If $\a\in C^0$, then  
$\G(\a)\geq -1$ which proves that $\a\in S$.\qed

%\bigskip
\vfill\eject

\centerline{\bf Replacing the Positive Solution by a Holomorphic Chain}

\Theorem{2.7} {\sl  Suppose $\G$ is a scarred $(2p-1)$-cycle (of class $C^1$) in an arbitrary complex manifold $X$.  Assume each component $\G_k$ of $\G$  can be pushed out at some point.
If \ $\G= dT$ for some weakly  positive $(p,p)$-current  $T$, then there exists a positive  holomorphic p-chain 
$\bbt$ with $\G=d\bbt$  and $T-\bbt\geq 0$,  so in particular,
$M(\bbt)\leq M(T)$ and $\supp(\bbt)\subset \supp(T)$.
}
\medskip

The proof depends on the following local result.

\smallskip
\Lemma {2.8}  {\sl Suppose $\G$ is an oriented connected $(2p-1)$-dimensional submanifold near $0\in \G$ in $\bbc^n$.   
\smallskip

\item
{(1)}  If $\G$ can be pushed out at $0\in \G$ and $r\G=dT$ for some $T\geq 0$ and $r>0$, then
$\G$ is two-sided near 0.  That is,  near 0 there exists a (unique) complex $p$-dimensional subvariety $V$ containing $\G$, so that $V=V^+\cup \G \cup V^-$ and $dV^{\pm} = \pm \G$.

\smallskip
\item
{(2)}  If $\G$ is two-sided near 0 and $r\G=dT$ for some $T\geq 0$ and $r>0$, then}

$$T=rV^+ +S \qquad {\rm with}\ \   S\geq 0 {\ \ and\ \ } dS=0.$$

\noindent
{\bf Proof.} By the push-out hypothesis  we have that $-\G=dZ$ for some irreducible subvariety $Z$ of $B(0, R)-\G$.
By taking a small piece $V^-$ of $Z$ we may assume that the positive current $T^+ \equiv T+rV^-\geq0$ has boundary $\G^+$ which does {\sl not} contain the origin. (See Figure 1.)
Consider the subset
$$
E_r(T^+) \ =\ \{z: \Theta (T^+, z)\geq r\}\subset B(0,R)-\G^+
$$
where $ \Theta (T^+, z)$ denotes the standard density, or Lelong number, of $T^+$ at $z$.
Since $dT^+=0$ in $B(0,R)-\G^+$ we know by a fundamental theorem of Siu [Siu] that
\medskip

\centerline{{\sl $E_r(T^+) $ is a complex subvariety of complex dimension $\leq p$
and }}\medskip

\centerline{{\sl $T^+-rW\geq0$ \ \  where $W$ is the $p$-dimensional part of $E_r(T^+) $.}}
\medskip

\noindent
  Since $E_r(T^+)$ contains
$V^-$, it must have an irreducible $p$-dimensional component $V\supset V^-$, defined 
in a neighborhood of the origin. This proves (1).  Since $V\subset W$, we have $T^+-rV\geq 0$.
Note also that $d(T^+ -rV)=0$ near the origin.  This proves (2) since
$S\equiv  T-rV^+ = T^+- rV$.\qed

\Cor{2.9} {Under the hypotheses of Lemma 2.8 (1), $-\G$ can also be pushed out at 0.}

\medskip\noindent
{\bf Proof of Theorem 2.7.}  As an easy consequence of Siu's Theorem (See, for example, Theorem 2.4, p. 638 in [H$_1$]), 
there exist irreducible $p$-dimensional subvarieties $V_j$ of $X-\supp \G$ and positive constants $c_j$ so that 
$$
T\ =\ \sum_{j=1}^{\infty} c_j V_j +R
\eqno{(2.3)}
$$
where $R\geq 0$ and, for each $c>0$, the complex subvariety $E_c(R)$ is of dimension $\leq p-1$.
This representation (2.3) is unique.  (Note that $R\geq 0$ implies that the mass of $T$ dominates the mass of $\sum_j c_j V_j$ on any set.)

Near the point where $\G_1$ can be pushed out, Lemma 2.8 (with $r=n_1$) implies that
$$
T\ =\ n_1V^+ + S \qquad\  {\rm with\ \ } S\geq 0 \ \ {\rm and \ \ }dS=0.
\eqno{(2.4)}
$$
By uniqueness $V^+$ must be contained in one of the $V_j$, say $V_1$.  Moreover, 
since $S\geq 0$ we have $c_1\geq n_1$.  This implies $\wt T \equiv T-n_1 V_1\geq 0$.

Near the point where $\G_1$ can be pushed out we have $dV_1=\G_1$. Hence, on $X$ we have
$dV_1 = \G_1 +\sum_{k=2}^{\infty} m_k \G_k$  with $m_k\in\bbz$.  Consequently, 
 $$
 d\wt T = \sum_{k=2}^{\infty}  (n_k - n_1 m_k)  \G_k,
 $$
 and so we have eliminated one of components of the boundary.
Now   the coefficients  in this sum may not all be positive, and to make them
all positive we may have to reverse the orientation of some of the $\G_k$.
However, by Corollary 2.9 these orientation-reversed components can also
be pushed out at some point. Hence, $\wt \G = d\wt T = \sum_{k=2}^{\infty} {\wt n}_k \wt {\G_k}$ 
satisfies all the hypotheses of Theorem 2.7.

If  $\G$ has only a finite number of components, then we are done by induction on the number 
of components.  If not, then by continuing this process we obtain a sequence of positive currents
${\wt T}_k = T-(n_1V_1+n_2'V_2+\cdots +n_k'V_k)$  where the $n_j'$ are positive integers and
$$
d {\wt T}_k \ =\ \sum_{j=k+1}^\infty n_{kj}\G_j.
$$
Since $T-{\wt T}_k  \geq 0$, we may assume, by passing to a subsequence, that $\{ {\wt T}_k\}_{k=1}^\infty$ converges in mass norm to a positive  current $ {\wt T}_\infty$, which must be flat since each
$ {\wt T}_k$  is a normal current.  Note that  $\supp(d {\wt T}_\infty)\subset \Scar$  and recall that,
by assumption, the scar set $\Scar$ has    Hausdorff (2p-1)-measure zero.   Hence, by 
%the  Federer Flat Support Theorem 
[F, 4.1.20],  we have $d {\wt T}_\infty=0$.
We conclude that $V = \sum n_j' V_j = T-{\wt T}_\infty$ is a positive holomorphic chain with 
$dV = \G$.
 \qed

\medskip
\noindent
{\bf Proof of Theorems  1.1 and 1.9.}
Remarks 1.5 and 2.5 show that if $\G=\sum_k m_k\G_k$ satisfies the linking hypothesis, then $\G_k$ is two-sided at any real analytic point.

\medskip
\noindent
{\bf Proof of Theorem 1.6.}
It suffices to show that when $X$ is $q$-convex for $q\leq p$, then Theorem 2.6 also holds 
with $\G$ and $T$ having compact support. For this we change the definitions of $S$  and $C$
in the proof of Theorem 2.6 by permitting the $\a$'s in $S$ to have arbitrary support and 
restricting the $T$'s in $C$ to have compact support.
 The argument will  carry through as before 
once it is established that the cone $C$  is closed in the weak topology.
This follows from standard compactness theorems and the following fact.  Suppose $f:X\to\bbr^+$
is the proper exhaustion  with $n-q+1$ positive eigenvalues on $\{x : f(x)\geq 1\}$. 
If  $T\in \cp_{p,p}(X)$  for $p\geq q$, then
$$
\supp\, T \ \subset  \biggl\{x\in X: f(x)\leq \max\biggl\{1, \sup_{\supp\, dT} f \biggr\}  \biggr\}.
$$

\vfill\eject

%%%%%%%%%%%%%%%%%%%%%%%%%%%%%%%%%%%%%%%%%%%%
%%%%%%%%%%%%%%%%%%%%%%%%%%%%%%%%%%%%%%%%%%%%
%%%%%%%%%%%%%%%%%%%%%%%%%%%%%%%%%%%%%%%%%%%%
%%%%%%%%%%%%%%%%%%%%%%%%%%%%%%%%%%%%%%%%%%%%
%%%%%%%%%%%%%%%%%%%%%%%%%%%%%%%%%%%%%%%%%%%%

\noindent{\bf \CC. Relative Hodge Classes and Representability.} \smallskip

In this chapter we address the question of  when a relative homology class 
can be represented by  a positive holomorphic chain.  
More specifically,  let $X$ be a compact K\"ahler manifold and $M\ss X$
a smooth orientable compact  submanifold  of real dimension $2p-1$.
Then we have the two closely related questions:

\medskip
\noindent
{\bf Relative Hodge Question:} {\sl  Which classes in $H_{2p}(X,M;\,\bbz)$ can be represented
by holomorphic chains?}

\medskip
\noindent
{\bf Relative Hodge Question (Positive Version):} {\sl  Which classes in $H_{2p}(X,M;\,\bbz)$ can be represented by  positive holomorphic chains?}

\medskip
We shall work in the space $\overline{H}_{2p}(X,M;\,\bbz) = H_{2p}(X,M;\,\bbz)/{\rm Tor}$  
where Tor is the torsion subgroup and  use
the following Relative de Rham Theorem.
Consider the short exact sequence of chain complexes of Fr\'echet spaces
$$
0\ \arr\  \ce^*(X,M)\ \arr\ \ce^*(X)\ \harr {j^*} \ \ \ce^*(M)\ \arr\ 0,
\eqno{(\CC.1)}
$$
where $j:M\to X$ denotes the inclusion, 
and the dual sequence of topological dual spaces
$$
0\ \longleftarrow\  {\ce_*'(X)\over  j_*\ce_*'(M)}\ \longleftarrow \ \ce_*'(X)\ \leftharr {j_*} \ \ \ce_*'(M)\ \longleftarrow\ 0.
\eqno{(\CC.2)}
$$
The complex $\ce^*(X,M)$, consisting of    forms which vanish when restricted to $M$, computes 
the relative cohomology $H^*(X,M;\,\bbr)$, and the complex 
$\ce_*'(X,M) \equiv \ce_*'(X)/  j_*\ce_*'(M)$ computes the relative homology
$H_*(X,M;\,\bbr)$.

The Relative de Rham Theorem states that:\medskip
\centerline
{\sl
$H^k(X,M;\,\bbr)$ and $H_k(X,M;\,\bbr)$ are dual to each other.}
\medskip
\noindent This can  be proven as follows.  Consider the dual triples
$$
\eqalign
{
 \ce^{k-1}(X,M) \  \harr  d \  \  &\ce^{k}(X,M) \  \harr d \ \  \ce^{k+1}(X,M)   \cr
 \ce_{k-1}'(X,M) \ \leftharr{d} \ \  &\ce_{k}'(X,M)\ \leftharr {d} \ \  \ce_{k+1}'(X,M)   \cr
}
$$
where $H^k(X,M;\,\bbr) = Z/B$ using the cycles $Z$ and boundaries $B$ in the first sequence,
and $H_k(X,M;\,\bbr) = \wt Z/\wt B$ using the cycles $\wt Z$ and boundaries $\wt B$ in the 
second sequence.  By the Hahn-Banach Theorem it suffices to show that $B$ and $\wt B$ are
closed.  These spaces are images of continuous linear maps.  If they  are of finite codimension
in $Z$ and $\wt Z$ respectively, then they are closed by a standard result in functional analysis.
Thus is remains to show that $H^k(X,M;\,\bbr)$ and $H_k(X,M;\,\bbr)$ are finite dimensional.
That $H^k(X,M;\,\bbr)$ is finite dimensional  follows  from  the long exact sequence
$$
\cdots \arr \ H^{k-1}(M;\,\bbr) \ \arr\ H^{k}(X,M;\,\bbr) \ \arr\ H^{k}(X;\,\bbr)\ \arr\ H^{k}(M;\,\bbr)\ \arr\ \cdots
$$
derived from (\CC.1) and the fact that  $H^{k-1}(M;\,\bbr)$ and  $H^k(X;\,\bbr)$ are finite dimensional
by the standard de Rham Theorem. That $H_k(X,M;\,\bbr)$ is finite dimensional follows similarly
from the long exact sequence derived from (\CC.2).

In the special case $k=2p = \dim M+1$ we have:
$$
\eqalign
{H^{2p}(X,M;\,\bbr) &= Z/B \qquad {\rm where}  \cr
Z\ =\ \{\vf \in\ce^{2p}(X) : d\vf=0\} \quad &{\rm and}\quad B\ =\ d\ce^{2p-1}(X,M)
}
\eqno{(\CC.3)}
$$
and $$
\eqalign
{H_{2p}(X,M;\,\bbr) = \wt Z/\wt B \qquad &{\rm where}  \cr
\wt Z\ =\ \{ T  \in\ce_{2p}'(X) : d T \in j_*\ce_{2p-1}'(M)\} \quad &{\rm and}\quad
 \wt B\ =\ d\ce_{2p+1}'(X).
}
\eqno{(\CC.4)}
$$

It is an interesting fact, established in the next section, that the group 
$\overline{H}_{2p}(X,M;\,\bbz) \ss  {H}_{2p}(X,M;\,\bbr)$ carries a ``real Hodge filtration''.
A key point is the following lemma. 

\Lemma{\CC.1}  {\sl  Fix $\tau \in  {H}_{2p}(X, M;\,\bbr)$.
If $T, T'\in \ce_{2p}(X,M)$ are relatively closed currents  representing $\tau$, then
$$
T(\vf) \ =\  T'( \vf)   \quad {\sl for\ all}\ \  \vf \in \ce^{2p}(X)\ \ {\sl with}\ \ d\vf=0.
\eqno{(\CC.5)}
$$
Hence, the notion of $\tau(\vf)$ is well defined for such $\vf$.
Furthermore,
$$
d T \ =\  d T'\ =\  \sum_{j=1}^L r_j [M_j]
\eqno{(\CC.6)}
$$
where  $M=M_1\cup\cdots \cup M_L$ is the decomposition into
connected components and the  $r_j$ are real numbers.  Thus,
$\partial \tau \ =\ \sum_{j=1}^L r_j [M_j]$ is well defined.
}

\pf  
Since $T$ and $T'$ both represent $\tau\in \wt Z/\wt B$ and 
$\wt B = d\ce_{2p+1}'(X)$ by (\CC.4), we have $T-T'= dR$
for $R\in \ce_{2p+1}'(X)$. This proves (\CC.5) and that $dT=dT'$.
Since $T\in \wt Z$, (\CC.4) says that $dT=j_*u$ where $u\in \ce_{2p-1}(M)$.
This implies that $du=0$.  Hence, $u$ is a locally constant function on $M$.
\qed
\medskip

\Def {\CC.2} A class $\tau \in  H_{2p}(X, M;\,\bbr)$  is called {\bf  positive} if $\tau(\vf)\geq 0$ for all closed, real 
$2p$-forms $\vf$ such that the component
$$
\vf^{p,p} \  \geq\ 0\quad {\rm (is\ weakly\ positive)\ on\ }  X.
$$
If $\tau$ is positive, then it is of type $(p,p)$ as defined in $\EE.1$ below.

\Prop{\CC.3} {\sl A class $\tau\in H_{2p}(X, M;\,\bbr)$ is positive if and only if it 
is represented (in the complex $\ce_*(X,M)$) by a strongly positive current of type $(p,p)$.}
\medskip

This proposition  will be proved below.   We first observe that it leads to the following
main result.

\Theorem {\CC.4}  {\sl Suppose $\tau \in  \overline{H}_{2p}(X, M;\,\bbz)$ is positive.
Suppose each component of $M$ has a real analytic point (or, more generally,
 is two-sided at some point).  
Then there exists a positive holomorphic $p$-chain $T$ on $X$
with $d T=\partial \tau$. Furthermore, there exists a positive
$d$-closed $(p,p)$-current $S$ with $\tau =[T+S]$.

In particular, if the positive  classes in $H_{2p}(X;\,\bbq)$ are all 
represented  by positive holomorphic chains with rational coefficients, then so are all the positive
classes in $H_{2p}(X, M;\,\bbq)$.}

\medskip

Thus for example,  given any real analytic $M$ in a Grassmann manifold $X$, we conclude that
every positive class in ${H}_{2p}(X, M;\,\bbz)$ carries a positive holomorphic chain.
However there are projective manifolds $X$ with positive $(p,p)$-classes in 
$\overline{H}_{2p}(X;\,\bbz)$ which do not carry
positive holomorphic cycles.  In fact, 
  for every integer $k\geq 2$ there exists an abelian variety $X$ of complex dimension $2k$
and a class $\tau \in H_{2k}(X;\,\bbz)$ which is represented by a positive $(k,k)$-current
and also by an algebraic $k$-cycle, but $\tau$ is not represented by a positive algebraic $k$-cycle
(see [L]).

\pf
By Proposition \CC.3 and (\CC.6) in Lemma \CC.1, 
the class $\tau$ is represented by a positive $(p,p)$-current  $T$
with $d T=\partial \tau = \sum_i n_i [M_i]$ for integers $n_i$ (cf. the argument for (\CC.2) above.)
Applying Theorem 2.7 with $\Gamma=d T$, we deduce the existence of 
a positive holomorphic chain $T_0$ with 
$$
\qquad\qquad\qquad\qquad\qquad
d T_0 \ =\ d T \and T-T_0\geq 0. 
\qquad\qquad\qquad\qquad\qquad 
\vrule width5pt height5pt depth0pt
$$ 
 
 \medskip\noindent
 {\bf Proof of Proposition \CC.3.}
Consider  the closed convex cones
$$\eqalign
{
P\ &\equiv\  \{\vf \in\ce^{2p}(X) : \vf ^{p,p} \ {\rm is\  weakly\ positive}\}\ \ss \ \ce^{2p}(X)  \cr
\wt P \ &\equiv\ \{T\in \ce_{2p}'(X) : T= T_{p,p}  \ {\rm is\  strongly\ positive}\}\ \ss\ \ce_{2p}'(X) \cr
}
$$
These are polars of each other in the dual pair $\ce^{2p}(X)$, $\ce_{2p}'(X)$.
Moreover, by the Relative de Rham Theorem and (\CC.3) and (\CC.4) we have:

\medskip

\item{(i)} \ \ $\wt B \ss\ce_{2p}'(X)$ is closed (in the weak topology).

\smallskip

\item{(ii)} \ \ $Z$ and $\wt B$ are polars of each other in the dual pair $\ce^{2p}(X)$, $\ce_{2p}'(X)$.

\smallskip

\item{(iii)} \ \ $B \ss\ce^{2p}(X)$ is closed.

\smallskip

\item{(iv)} \ \ $B$ and $\wt Z$ are polars of each other in the dual pair $\ce^{2p}(X)$, $\ce_{2p}'(X)$.

\Lemma {\CC.5} {\sl The subset    $\wt P+\wt B$ is closed in the standard
topology on  $\ce_{2p}'(X)$.}

\pf  Let $\{T_i\}\ss \cp_{p,p}$ and $\{d S_i\}\ss \wt B$ be sequences such that 
$$
T_i+d S_i\ \arr\ R\qquad{\rm weakly\ in\ \ } \ce_{2p}'(X)
$$
Let $\o$ denote the K\"ahler form on $X$.  Then
$$
p!M(T_i) \ =\ T_i\left( \omega^p\right) \ 
=\ (T_i+d S_i)\left( \omega^p\right) \ \arr\ R\left(  \omega^p\right),
$$
and so the masses $M(T_i)$ are uniformly bounded. By
the compactness theorem for positive currents there is a subsequence, 
again denoted by $T_i$, converging 
to a  positive current $T$. Hence, $d S_i \arr R-T$ weakly, and since
$d$ has closed range, there exists $S\in \ce_{2p+1}'(X)$ with $d S=R-T$.\qed

\Prop{\CC.6}
$$
\left[(P\cap Z)+B\right]^0 \ = \ (\wt P\cap \wt Z)+ \wt B
$$

\pf
By standard principles we have
$\left[(P\cap Z)+B\right]^0 = (P\cap Z)^0 \cap B^0 = 
\overline{(P^0+Z^0)} \cap B^0
$.
By (ii), (iv) and  Lemma \CC.5  we have
$\overline{(P^0+Z^0)} \cap B^0 = \overline{(\wt P+\wt B)} \cap \wt Z = {(\wt P+\wt B)} \cap \wt Z
$.
Finally it is easy to see that 
${(\wt P+\wt B)} \cap \wt Z= (\wt P\cap \wt Z) +\wt B$ since $\wt B\ss\wt Z$.\qed

\medskip
To complete the proof of Proposition \CC.3
choose a current $T\in \wt Z$ which represents the class $\tau$.  By hypothesis
$T$ is in the polar of  $(P\cap Z)+B$.  Therefore, by Proposition \CC.6 and (\CC.4),  $T=T_0+d S$
with  $T_0\in \wt P$.\qed

 %%%%%%%%%%%%%%%%%%%%%%%%%%%%%%%%%%%%%%%%%%%%%%
%%%%%%%%%%%%%%%%%%%%%%%%%%%%%%%%%%%%%%%%%%%%%%
%%%%%%%%%%%%%%%%%%%%%%%%%%%%%%%%%%%%%%%%%%%%%%
\vskip .3in

\noindent
{\bf \S \EE.  A Real Hodge Filtration on   $H_{2p}(X,M;\,\bbr)$.}

\Def{\EE.1}    A homology class $\tau
\in H_{2p}(X, M;\,\bbr)$ is of  {\bf  filtration level $k$}  if $\tau(\vf)= 0$  for
all closed complex valued forms $\vf$ of type $(r,s)$ with $r>p+k$.
Classes of filtration level 0 are called {\bf type (p,p)}.

\Note{\EE.2} This induces a real Hodge filtration $F^kH_{2p}(X,M;\,\bbr)$ on $H_{2p}(X, M;\,\bbr)$
which extends the basic one $F^kH_{2p}(X;\,\bbr) = 
\bigoplus_{r=0}^k \left\{ H_{p-r, p+r}(X)\oplus H_{p+r, p-r}(X)\right\}_\bbr$ on 
 $H_{2p}(X;\,\bbr)$.

\Prop{\EE.3}  {\sl Suppose $\tau \in H_{2p}(X,M;\,\bbr)$ has filtration level $k$. 
Then $\tau$ is represented by a current
$$
T\ \in \ \left\{ \ce_{p-k,p+k}'(X) \oplus  \cdots \oplus \ce_{p+k,p-k}'(X)
\right\}_{\bbr},
\eqno{(\EE.1)}
$$
and therefore,
$$
d T \ \in \ \left\{ \ce_{p-k-1,p+k}'(X) \oplus  \cdots \oplus \ce_{p+k,p-k-1}'(X)
\right\}_{\bbr}.
\eqno{(\EE.2)}
$$
In particular, if $\tau$ is of type $(p,p)$, then $\tau=[T]$ for a some $(p,p)$-current $T$,
and each non-zero boundary component of $\partial \tau$ is maximally complex
(cf. [HL$_1$]).
}

\pf  We start by establishing (\EE.2). 
 Write $\partial \tau = \sum_j r_j [M_j]$ as in Lemma \CC.1.  Choose any smooth form $\psi \in 
\ce^{r,s}(X)$ with $r+s=2p-1$ and either $r>p+k$ or $s>p+k$.   Then $0=\tau(d\psi) = (\partial \tau)(\psi)
= \sum_j r_j \int_{M_j} \psi$.  Since $\psi$ is arbitrary, we conclude that for each $M_j$ with 
$r_j\neq 0$, the Dolbeault components 
$$
[M_j]_{r,s}=0 \quad {\rm if \ either\ \ }  s>p+k\ \ {\rm or}\ \ r>p+k.
$$
This gives (\EE.2).  When $k=0$ this means  $M_j$ is maximally complex.

Consider the case where $\tau$ is of type $(p,p)$ with $2p\leq n$.
Choose a current $T$ representing $\tau$. Then by standard harmonic theory
$T_{2p,0} = h_{2p,0}-\dbar \b$ where $h$ is  harmonic (in particular, smooth)
and $\b\in \ce_{2p,1}'(X)$.  Then $[T-d\b ]_{2p,0} = T_{2p,0} - \dbar \b =  h_{2p,0}\equiv h$
and because $T(\bar{*}h) = \|h\|^2 =0$ (since $\tau=[T]$ is type $(p,p)$), we have $h=0$.
Thus replacing $T$ by $T-\dbar \b-\partial \overline \b$ we can assume $T_{2p,0}=T_{0,2p}=0$.

If $p=1$, we are done. If $p>1$, we note that $\dbar T_{2p-1,1} = M_{2p-1,0} = 0$, and so
$T_{2p-1,1} = h_{2p-1,1} + \dbar \b$ where $ h_{2p-1,1}$ is harmonic and $\b \in \ce_{2p-1,2}'(X)$.
We conclude as above that $ h_{2p-1,1} = 0$, and then replace $T$ by $T-\dbar \b-\partial \overline \b$
so that  $T_{2p-1,1}=T_{1,2p-1}=0$. Continuing in this fashion gives the result. All other cases
are entirely analogous and details are left to the reader.\qed

\vfill\eject

%\vskip .3in

\noindent
{\bf \S \DD.  Residual Currents.}
\Def{\DD.1} Let $R$ be a  weakly positive, $d$-closed $(p,p)$-current. Then  $R$ is {\bf residual } if for each $c>0$ the complex dimension of the subvariety $E_c(R) = \{z : \Theta(R,z)\geq c\}$ is $\leq p-1$.

\medskip
Suppose $T$
is a  a  weakly positive, $d$-closed $(p,p)$-current defined in the complement of supp$(\G)$
where $\G$ is a scarred $2p-1$ cycle (of class $C^1$).  By the main result of [H]
(see Theorem 6, p. 71 and the note added in proof)  $T$ has locally finite mass across 
supp$(\G)$.  That is, $T$ has a unique ``extension by zero'' across supp$(\G)$. Let $T$ also denote this extension.  It follows easily that from two theorems of Federer that $dT=\sum_{k=1}^\infty r_k \G_k$
with constants $r_k\in\bbr$.

\Def{\DD.2}   The set  supp$(\G)$ is {\bf residual} if each residual current $R$ on $X-\supp(\G)$ satisfies $dR=0$ on $X$.

\Prop{\DD.3}  {\sl If each component of $\G$ has a two-sided point, then $\supp(\G)$ is residual.}
\pf
Suppose $R$ is a residual current on $X-\supp(\G)$ with $dR=\sum_k r_k\G_k$, $r_k\in\bbr$. Near a two-sided point of one of the components, say$\G_1$, we have $dR=r_1\G_1$.
By Lemma 2.8 we can write $R$ locally as  $R=r_1V^++S$ with $S\geq 0$ and $dS=0$ across $\G_1$. This contradicts the hypothesis that $R$ is residual unless $r_1=0$. \qed

\medskip
\noindent
{\bf Remark.}  This Proposition combined with the first half of Lemma 2.8 and the next result provides a second proof of Theorem 2.4.

\Theorem{\DD.4} {\sl Suppose $\G$ is a scarred $2p-1$ cycle (of class $C^1$) in an arbitrary complex manifold $X$.  Assume each component $\G_k$ of $\G$ is residual at some point.  If $\G=dT$ for some weakly positive $(p,p)$-current $T$ on $X$, then there exists a positive holomorphic $p$-chain $\bbt$
with $\G=d\bbt$ and $T-V\geq0$.}

\pf  Suppose $\G=dT$ as in the theorem and consider the decomposition $T=S+R$ into a positive
{\sl real-coefficient holomorphic chain} $S=\sum_{j=1}^\infty c_j V_j$ plus a residual current $R$
(on $X-\G$),  Now $dR = \sum_{k=1}^\infty r_k \G_k$ for some $r_k\in \bbr$, but by the hypothesis
each $r_k$ must be zero.  Hence, $\G=dS$ bounds a positive real-coefficient holomorphic chain.

\Prop{\DD.5}  {\sl Let $\G$ be a scarred $2p-1$ cycle in an arbitrary complex manifold $X$.
If $\G=dS$ bounds a positive real-coefficient holomorphic chain $S=\sum_{j=1}^\infty c_j V_j$,
then  $\G=dV$ bounds a positive (integer-coefficient) holomorphic chain $V$ with $S-V\geq0$ 
(and therefore also  $\supp V\subseteq \supp S$).}

\pf By hypothesis $d\left(\sum_{j=1}^\infty c_j V_j\right) = \sum_{k=1}^\infty n_k \G_k$. Near a regular point $x$ in $\G_1$ each $V_j$ satisfies $dV_j = \epsilon_j\G_1$ with $\epsilon \in \{-1,0,1\}$.  By uniqueness there is at 
most one of the subvarieties $V_j$ with boundary $\G_1$.  Relabel so that $dV_1=\G_1$.
Now there are two cases.
\smallskip
\noindent
{\bf Case 1:}   $dV_j=0$  for all $j\geq 2$.
In this case we must have $c_1=n_1$, and we can eliminate the component $\G_1$ from $\G$.

\smallskip
\noindent
{\bf Case 2:}   $-\G_1$  bounds exactly one of the subvarieties $V_j$, $j\geq 2$.
Relabel so that $-\G_1 = dV_2$.  In this case $c_1-c_2=n_1$.  Note that $V_1+V_2$
is a subvariety without boundary near the point $x$ on $\G_1$.  Set $\wt S= S-n_1V_1 = c_2(V_1+V_2) + \sum_{j=3 }^\infty  c_j V_j$.  Then $\wt S$ is positive and $d\wt S=0$ near the point $x$.  Consequently, $d\wt S = \sum_{j=2 }^\infty  b_j V_j$.  Finally, the $b_j$'s must be integers.  In fact
$b_j=n_j-\epsilon_{1j }n_1$  where $dV_1=\epsilon_{1j}\G_j$ defines $\epsilon_{1j}\in\{-1,0,1\}$.
Hence,  we can eliminate the component $\G_1$ from $\G$ in this case as well
\smallskip

The proof can now  be completed exactly as in the last paragraph of the proof of Theorem 2.7.
\qed

\Qu{\DD.6}  Which (maximally complex) $(2p-1)$-dimensional submanifolds are residual?
Note that if $\G$ is two-sided, then $\G$ is residual.  Moreover, if $\G$ is one-sided, then $\G$ has a natural orientation so that $\G=d W$ where $W$ is complex, and in this case the residual property is equivalent to the following uniqueness property:
\medskip
\centerline{If $T\geq 0$ satisfies $dT=\G$, then $T=W+S$ with $S\geq0$ and $dS=0$.}
\medskip\noindent
If $\G$ is zero-sided, then $\G$ is residual if and only if 
$$
T\geq 0\ \ \ {\rm and }\ \ \ \supp\{dT\}\ \subset\ \G \qquad\Rightarrow\qquad dT=0.
$$

\vskip .4in

 % \magnification=1200
%\NoBlackBoxes
%\nologo 

%  \input  Defs.tex

\centerline{\bf References}

\vskip .1in

\noindent
\item{[AW]}  H. Alexander and J. Wermer, {\sl Linking numbers
and boundaries of varieties}, Ann. of Math.
{\bf 151} (2000),   125-150.

 \smallskip

\noindent
\item{[D]}  J.-P. Demailly, {\sl Estimations $L^2$ pour
l'op\'erateur $\overline{\partial}$ d'un fibr\'e vectoriel holomorphe
semi-positif au-dessus d'une vari\'et\'e k\"ahl\'erienne compl\`ete},
   Ann. Sci. E. N. S.  
 {\bf 15} no. 4  (1982),  457-511.

 \smallskip

 \noindent
\item{[Do]}  P. Dolbeault,   {\sl On holomorphic chains with given 
boundary in $\bbc\bbp^n$},    Springer  Lecture Notes in Math., {\bf 
1089}  (1983), 1135-1140.

\smallskip

\noindent
\item{[DH$_1$]}   P. Dolbeault and G. Henkin ,  
{\sl  Surfaces de Riemann de bord donn\'e dans $\bbc\bbp^n$  },  pp.
163-187 in  `` Contributions to Complex Analysis and Analytic Geometry''  ,
  Aspects of Math.  Vieweg {\bf 26 },   1994.

\smallskip

\noindent
\item{[DH$_2$]}   P. Dolbeault and G. Henkin,   
{\sl Cha\^ines holomorphes de bord donn\'e dans $\bbc\bbp^n$},    
Bull.  Soc. Math. de France,
{\bf  125}  (1997), 383-445.

\smallskip

\noindent
\item{[F]}   H. Federer, Geometric Measure  Theory,
 Springer--Verlag, New York, 1969.

 \smallskip

\noindent
\item{[GZ]}  V. Guedj and A. Zeriahi,     
{\sl    Intrinsic capacities on compact K\"ahler manifolds},    
Preprint Univ. de Toulouse , 2003

\smallskip

\noindent
\item{[H$_1$]}  F.R. Harvey,
{\sl Three structure theorems in several complex
variables},  Bull. A. M. S. {\bf 80} (1974), 633-641.

 \smallskip

\noindent
\item{[H$_2$]}  F.R. Harvey,{\sl
Holomorphic chains and their boundaries}, pp. 309-382 in ``Several Complex
Variables, Proc. of Symposia in Pure Math. XXX Part 1'', 
A.M.S., Prov., RI, 1977.

 \smallskip

\noindent
\item{[HK]}  F.R. Harvey and A. Knapp  {\sl Positive $(p,p)$-forms, Wirtinger's inequality, and
currents}, Part A (Proc. Tulane University Program on Value-Distribution Theory in Complex 
Analysis and Related Topics in Differential Geometry, 1972-73), pp. 43-62.  Dekker, New York, 1974.

 \smallskip
\noindent
\item{[HL$_1$]} F. R. Harvey and H. B. Lawson, Jr, {\sl On boundaries of complex
analytic varieties, I}, Annals of Mathematics {\bf 102} (1975),  223-290.

 \smallskip

 \noindent 
\noindent
\item{[HL$_2$]}  F. R. Harvey and H. B. Lawson, Jr, {\sl On boundaries of complex
analytic varieties, II}, Annals of Mathematics {\bf 106} (1977), 
213-238.

 \smallskip

 \noindent
\item{[HL$_3$]} F. R. Harvey and H. B. Lawson, Jr, {\sl Boundaries of 
varieties in projective manifolds}, J. Geom. Analysis,  {\bf 14}
no. 4 (2005), 673-695. ArXiv:math.CV/0512490
 \smallskip

 \noindent 
\noindent
\item{[HL$_4$]}  F. R. Harvey and H. B. Lawson, Jr, {\sl Projective hulls and
the projective Gelfand transformation}, Asian J. Math. (to appear).
ArXiv:math.CV/0510286

 \smallskip

\noindent
\item{[HL$_5$]} F. R. Harvey and H. B. Lawson, Jr, {\sl Projective  linking and 
boundaries of positive holomorphic chains in 
projective manifolds, Part I}, ArXiv:math.CV/0512379.

 \smallskip

\noindent
\item{[HLZ]} F. R. Harvey, H. B. Lawson, Jr. and J. Zweck, {\sl A
deRham-Federer theory of differential characters and character duality},
Amer. J. of Math.  {\bf 125} (2003), 791-847.  ArXiv:math.DG/0512251

 \smallskip

\noindent
\item{[L]}   H. B. Lawson, Jr, {\sl  The stable homology of a flat torus}, Math. Scand. 
{\bf 36} (1975), 49-73.

 \smallskip

   \noindent
\item{[S]}  H. H. Schaefer,  Topological Vector Spaces,    Springer Verlag,
New York,  1999.

\smallskip

   \noindent
\item{[Sh]}    B. Shiffman,    {\sl  On the removal of singularities of analytic sets},    Michigan Math. J.,
 {\bf 15}  (1968), 111-120.

\smallskip

   \noindent
\item{[Siu]}    Y.-T.    Siu,    {\sl  Analyticity of sets associated to Lelong numbers and the extension on closed positive currents},    Inventiones Math.,
 {\bf  27}  (1974), 53-156.

\smallskip

 \noindent
\item{[W$_1$]}   J. Wermer  {\sl    The hull of a curve in $\bbc^n$},    
Ann. of Math., {\bf  68}  (1958), 550-561.

\smallskip

 \noindent
\item{[W$_2$]}   J. Wermer    {\sl    The argument principle and
boundaries of analytic varieties},     Operator Theory: Advances and 
Applications, {\bf  127}  (2001), 639-659.

\smallskip

%%%%%%%%%%%%%%%%%%%%%%%%%%%%%%%%%%%%%%%%%%%%%%
%%%%%%%%%%%%%%%%%%%%%%%%%%%%%%%%%%%%%%%%%%%%%%
\end

\end